\documentclass{article}

\input{diagrams}

\usepackage{amstex}
\usepackage{latexsym}

\newcommand{\mcm}[3]{\newcommand{#1}[#2]{{\ensuremath{#3}}}}

\mcm{\mc}{1}{\mathcal{#1}}
\mcm{\mr}{1}{\mathrm{#1}}
\mcm{\mi}{1}{\mathit{#1}}
\mcm{\mb}{1}{\mathbf{#1}}
\mcm{\cat}{1}{\mc{#1}}
\mcm{\scat}{1}{\Bbb{#1}}
\mcm{\fcat}{1}{\mb{#1}} 
\mcm{\twid}{1}{\widetilde{#1}}
\newcommand{\url}[1]{\mbox{\tt #1}}

\mcm{\pr}{2}{(#1,#2)}
\mcm{\triple}{3}{(#1,#2,#3)}
\mcm{\range}{2}{#1,\,\ldots\,,#2}
\mcm{\tuplebts}{1}{( #1 )}
\mcm{\tuple}{3}{\tuplebts{\range{#1,#2}{#3}}}
\mcm{\bftuple}{2}{\tuplebts{\range{#1}{#2}}}
\mcm{\rttuple}{1}{\tuplebts{\,\ldots\,,#1}}
\mcm{\atuplebts}{1}{\langle #1 \rangle}
\mcm{\atuple}{3}{\atuplebts{\range{#1,#2}{#3}}}
\mcm{\abftuple}{2}{\atuplebts{\range{#1}{#2}}}
\mcm{\arttuple}{1}{\atuplebts{\,\ldots\,,#1}}
\mcm{\elt}{0}{\in}
\mcm{\sub}{0}{\subseteq}
\mcm{\without}{0}{\setminus}
\mcm{\such}{0}{\:|\:}
\mcm{\nat}{0}{\Bbb{N}}	
\mcm{\iso}{0}{\cong}
\mcm{\eqv}{0}{\simeq}
\mcm{\id}{0}{\mi{id}}
\mcm{\emptybk}{0}{\:\:}
\mcm{\blank}{0}{(\:\:)}
\mcm{\op}{0}{\mr{op}}
\mcm{\ftrcat}{2}{[#1,#2]}
\mcm{\Set}{0}{\fcat{Set}}
\mcm{\Sets}{0}{\fcat{Set}}
\mcm{\Cat}{0}{\fcat{Cat}}
\mcm{\Bicat}{0}{\fcat{Bicat}}
\mcm{\Hom}{0}{\mr{Hom}}
\mcm{\End}{0}{\mr{End}}
\mcm{\homset}{3}{#1(#2,#3)}
\mcm{\mtihom}{4}{#1(\range{#2}{#3};#4)} 
\mcm{\mtiendo}{2}{#1(#2;#2)}
\mcm{\dom}{0}{\mr{dom}}
\mcm{\cod}{0}{\mr{cod}}
\mcm{\Mon}{0}{\mb{Mon}}
\newcommand{\pf}{\textbf{Proof}}
\newcommand{\done}{\hfill\ensuremath{\Box}}

\mcm{\go}{0}{\rTo}
\mcm{\og}{0}{\lTo}
\mcm{\goby}{1}{\rTo^{#1}}
\mcm{\ogby}{1}{\lTo^{#1}}
\mcm{\goesto}{0}{\,\longmapsto\,}
\mcm{\slob}{3}{(#1 \goby{#2} #3)}
\mcm{\bktdslob}{3}{\slob{#1}{#2}{#3}}
\mcm{\vslob}{3}
	{\left.
	\begin{diagram}[height=1.5em]
	#1		\\
	\dTo>{#2}	\\
	#3		\\
	\end{diagram}
	\right.}
\mcm{\bktdvslob}{3}{\left( \vslob{#1}{#2}{#3} \right)}
\newarrow{Get}....>
\newarrow{Equals}=====
\mcm{\goiso}{0}{\goby{\diso}}
\mcm{\gobyciso}{3}{{#1}\from {#2} \goiso {#3}}
\mcm{\parpair}{2}{\stackrel{\rTo^{#1}}{\rTo_{#2}}}

\newenvironment{slopeydiag}
	{\begin{diagram}[size=2em]}
	{\end{diagram}}

\newenvironment{ntdiag}
	{\begin{diagram}[size=1.5em,noPS]}
	{\end{diagram}}
\mcm{\diso}{0}{\sim}
\mcm{\vdiso}{0}{\wr}

\mcm{\node}{0}{\bullet}
\mcm{\enode}{0}{\circ}
\mcm{\nl}{1}{\stackrel{#1}{\node}}

\mcm{\nll}{1}{\raisebox{-1ex}[0ex][0ex]{\ensuremath{\stackrel{\node}{#1}}}} 

\newarrow{Edge}---->

\mcm{\bee}{0}{\cat{B}}
\mcm{\beep}{0}{\cat{B'}}
\mcm{\ess}{0}{\cat{S}}
\mcm{\tee}{0}{\cat{T}}
\mcm{\stee}{0}{\scat{T}}
\mcm{\teep}{0}{\cat{T'}}
\mcm{\steep}{0}{\scat{T'}}
\mcm{\wun}{0}{\fcat{1}}

\mcm{\trip}{2}{\homset{\bee}{#1}{#2}}
\mcm{\tripp}{2}{\homset{\beep}{#1}{#2}}
\mcm{\wee}{1}{\scriptstyle{#1}}
\mcm{\rep}{1}{\trip{\dashbk}{#1}}
\mcm{\ftrbi}{2}{\ftrcat{#1}{#2}}
\mcm{\Pshf}{0}{\ftrbi{\bee^{\op}}{\Cat}}

\mcm{\blob}{0}{\scriptscriptstyle{\bullet}}
\mcm{\ust}{0}{{}^{*}}
\mcm{\ubl}{0}{{}^{\blob}}
\mcm{\stbk}{0}{\blank^{*}}
\mcm{\blbk}{0}{\blank^{\blob}}
\mcm{\Mnd}{0}{\triple{\stbk}{\eta}{\mu}}
\mcm{\Cartpr}{0}{\pr{\cat{S}}{\ust}}
\mcm{\spn}{3}{#2 \og #1 \go #3}
\mcm{\spaan}{5}{#2 \ogby{#4} #1 \goby{#5} #3}
\mcm{\gph}{2}{\spn{#1}{#2^{*}}{#2}}
\mcm{\graph}{4}{\spaan{#1}{#2^{*}}{#2}{#3}{#4}}
\mcm{\comp}{0}{\mi{comp}}
\mcm{\ids}{0}{\mi{ids}}
\mcm{\Slice}{0}{\cat{S}/C_0}
\mcm{\unit}{0}{\mi{unit}}
\mcm{\mult}{0}{\mi{mult}}
\mcm{\Imnd}{0}{\triple{\blbk}{\unit}{\mult}}
\mcm{\Icartpr}{0}{\pr{\Slice}{\ubl}}
\mcm{\Alg}{0}{\mb{Alg}}
\mcm{\Multicat}{0}{\fcat{Multicat}}
\mcm{\Multicatof}{1}{#1\!\!-\!\!\Multicat}
\mcm{\Graph}{0}{\fcat{Graph}}
\mcm{\Graphof}{1}{#1\!\!-\!\!\Graph}
\mcm{\ush}{0}{{}^{\sharp}}

\mcm{\bdry}{0}{\partial}
\mcm{\Tr}{0}{\fcat{Tr}}
\mcm{\Cpn}{1}{\pr{\Set/S_{#1}}{T_{#1}}}
\mcm{\Operadof}{1}{#1\!\!-\!\!\fcat{Operad}}
\mcm{\Strucof}{1}{#1\!\!-\!\!\fcat{Struc}}
\mcm{\Catof}{1}{#1\!\!-\!\!\Cat}
\mcm{\CatOF}{1}{\Catof{\mb{#1}}}
\mcm{\BicatOF}{1}{\mb{#1}\!\!-\!\!\Bicat}
\mcm{\PD}{1}{\fcat{PD}_{#1}}
\mcm{\pd}{1}{\fcat{pd}_{#1}}
\mcm{\TR}{1}{\fcat{TR}(#1)}
\mcm{\Gy}{0}{\fcat{Gy}}
\mcm{\Ssq}{0}{\fcat{Ssq}}
\mcm{\tooc}{0}{\fcat{2\!\!-\!\!Cat}}
\mcm{\Cay}{0}{\mb{Cay}}
\mcm{\thomset}{2}{\homset{\tee}{#1}{#2}}
\mcm{\cayset}{1}{\coprod\thomset{B}{#1}}
\mcm{\dcayset}{1}{\coprod_{B\elt\stee_0}\thomset{B}{#1}}
\mcm{\bkpr}{1}{{}^{(#1)}}

\mcm{\Sym}{0}{\mr{Sym}}

\newenvironment{morphs}{\begin{array}{rccc}}{\end{array}}
\newcommand{\lump}[1]{\subsubsection*{#1}}
\newcommand{\deflump}[1]{{\raggedright\vspace{1em}\textsc{#1}\nopagebreak\\ 
	\vspace{0em} \nopagebreak}}
\mcm{\oppair}{2}{\stackrel{\rTo^{#1}}{\lTo_{#2}}}
\mcm{\diagspace}{0}{\mbox{\hspace{2em}}}

\mcm{\ob}{0}{\mr{ob}\,}
\mcm{\dashbk}{0}{\mbox{---}}
\newcommand{\from}{:}
\mcm{\gobyc}{3}{{#1}\from {#2} \go {#3}}	
\mcm{\of}{0}{\raisebox{0.2mm}{\ensuremath{\scriptstyle\circ}}}

\newlength{\globwd}	
\newlength{\globht}	
\newlength{\globdh}	
\newlength{\tempmeas}	
\newlength{\globdrop}	
\newlength{\offset}	

\newcommand{\cell}[4]{\put(#1,#2){\makebox(0,0)[#3]{\ensuremath{#4}}}}

\newcommand{\initdims}[1]{
	\setlength{\unitlength}{1em}
	\setlength{\offset}{.5\unitlength}
	\settowidth{\globwd}{#1}
	\settoheight{\tempmeas}{#1}
	\addtolength{\tempmeas}{.4\unitlength}	
	\setlength{\globht}{.5\tempmeas}
	\setlength{\globdh}{.5\tempmeas}
	\addtolength{\globht}{1\offset}
	\addtolength{\globdh}{-1\offset}
	\setlength{\globdrop}{1\globdh}}
\newcommand{\sidespic}[1]{
	\settowidth{\tempmeas}{\ensuremath{#1}}
	\addtolength{\globwd}{\tempmeas}}
\newcommand{\abovepic}[1]{
	\settoheight{\tempmeas}{\ensuremath{#1}}
	\addtolength{\globht}{\tempmeas}
	\settodepth{\tempmeas}{\ensuremath{#1}}
	\addtolength{\globht}{\tempmeas}}
\newcommand{\belowpic}[1]{
	\settoheight{\tempmeas}{\ensuremath{#1}}
	\addtolength{\globdh}{\tempmeas}
	\settodepth{\tempmeas}{\ensuremath{#1}}
	\addtolength{\globdh}{\tempmeas}}
\newcommand{\present}[1]{
	\raisebox{-1\globdrop}[1\globht][1\globdh]{\makebox[1\globwd]{#1}}}

\newcommand{\pretwc}[5]
{\begin{picture}(4.2,3.4)(-0.1,-0.2)
\cell{-0.1}{1.5}{r}{#1}
\cell{4.1}{1.5}{l}{#2}
\cell{2}{3.2}{b}{#3}
\cell{2}{-0.2}{t}{#4}
\cell{2.2}{1.5}{l}{#5}
\qbezier(0,2)(2,4)(4,2)
\qbezier(0,1)(2,-1)(4,1)
\put(4,2){\vector(1,-1){0}}
\put(4,1){\vector(1,1){0}}
\put(2,2.5){\vector(0,-1){2}}
\end{picture}}

\mcm{\twc}{5}{
\initdims{\pretwc{}{}{}{}{}}
\sidespic{#1}
\sidespic{#2}
\abovepic{#3}
\belowpic{#4}
\present{\pretwc{#1}{#2}{#3}{#4}{#5}}}

\newcommand{\pretwcop}[5]
{\begin{picture}(4.2,3.4)(-0.1,-0.2)
\cell{-0.1}{1.5}{r}{#1}
\cell{4.1}{1.5}{l}{#2}
\cell{2}{3.2}{b}{#3}
\cell{2}{-0.2}{t}{#4}
\cell{2.2}{1.5}{l}{#5}
\qbezier(0,2)(2,4)(4,2)
\qbezier(0,1)(2,-1)(4,1)
\put(0,2){\vector(-1,-1){0}}
\put(0,1){\vector(-1,1){0}}
\put(2,2.5){\vector(0,-1){2}}
\end{picture}}

\mcm{\twcop}{5}{
\initdims{\pretwcop{}{}{}{}{}}
\sidespic{#1}
\sidespic{#2}
\abovepic{#3}
\belowpic{#4}
\present{\pretwcop{#1}{#2}{#3}{#4}{#5}}}

\newcommand{\pretwcv}[5]
{\begin{picture}(3.4,4.2)(0.8,0.9)
\cell{2.5}{5.1}{b}{#1}
\cell{2.5}{0.9}{t}{#2}
\cell{0.8}{3}{r}{#3}
\cell{4.2}{3}{l}{#4}
\cell{2.5}{3.2}{b}{#5}
\qbezier(2,5)(0,3)(2,1)
\qbezier(3,5)(5,3)(3,1)
\put(2,1){\vector(1,-1){0}}
\put(3,1){\vector(-1,-1){0}}
\put(1.5,3){\vector(1,0){2}}
\end{picture}}

\mcm{\twcv}{5}{
\initdims{\pretwcv{}{}{}{}{}}
\abovepic{#1}
\belowpic{#2}
\sidespic{#3}
\sidespic{#4}
\present{\pretwcv{#1}{#2}{#3}{#4}{#5}}}

\newcommand{\prethcv}[7]
{\begin{picture}(5,8.2)(0.5,-1.6)
\cell{3}{6.6}{b}{#1}
\cell{3}{-1.6}{t}{#2}
\cell{0.5}{2.5}{r}{#3}
\cell{5.5}{2.5}{l}{#4}
\cell{3}{4.2}{b}{#5}
\cell{3}{0.8}{t}{#6}
\cell{3.2}{2.5}{l}{#7}
\qbezier(3.5,6.5)(7,2.5)(3.5,-1.5)
\qbezier(2.5,6.5)(-1,2.5)(2.5,-1.5)
\put(2.5,-1.5){\vector(1,-1){0}}
\put(3.5,-1.5){\vector(-1,-1){0}}
\qbezier(1,3)(3,5)(5,3)
\qbezier(1,2)(3,0)(5,2)
\put(5,3){\vector(1,-1){0}}
\put(5,2){\vector(1,1){0}}
\put(3,3.5){\vector(0,-1){2}}
\end{picture}}

\mcm{\thcv}{7}{
\initdims{\prethcv{}{}{}{}{}{}{}}
\abovepic{#1}
\belowpic{#2}
\sidespic{#3}
\sidespic{#4}
\present{\prethcv{#1}{#2}{#3}{#4}{#5}{#6}{#7}}}

\begin{document}

\title{Basic Bicategories}
\author{Tom Leinster\\ \\
	\normalsize{Department of Pure Mathematics, University of
	Cambridge}\\ 
	\normalsize{Email: leinster@@dpmms.cam.ac.uk}\\
	\normalsize{Web: http://www.dpmms.cam.ac.uk/$\sim$leinster}}
\date{1 May, 1998\\
	\small Minor amendments September 1998}

\begin{titlepage}
\maketitle
\thispagestyle{empty}

\begin{abstract}
A concise guide to very basic bicategory theory, from the definition
of a bicategory to the coherence theorem.
\end{abstract}

\vfill
\tableofcontents
\end{titlepage}


\section{Introduction}

This is a minimalist account of the coherence theorem for bicategories. The
definitions of a bicategory, of a morphism between them, and so on, are given
first; from here, a straight-line path is taken to the coherence theorem. No
motivation or context is given, and only such examples as are necessary to the
development of the theory. More discursive literature on bicategories is
available: for instance, the original paper of B\'{e}nabou, or Gray's book
(see bibliography). In particular, section~9 of Street's paper~\cite{StCS}
covers much the same material as this paper. 

Nothing here is new (although I haven't seen~\ref{sec:coh-and-comm} put this
way before). The definitions are culled from B\'{e}nabou's paper and Gray's
book, and the bare bones of the coherence theorem from the papers of
Street~\cite{StFB} and of Gordon, Power and Street. Many points are also
covered in Lack's thesis.

\subsection*{Acknowledgements}

This work was supported by a PhD grant from EPSRC. It was prepared in \LaTeX\
using Paul Taylor's diagrams package. I am grateful to Jeff Egger for
initiating the series of talks which induced me to learn about bicategories.
 


\section{Definitions}	\label{sec:defns}
\setcounter{subsection}{-1}

\subsection{Bicategories}

A \emph{bicategory} \bee\ consists of the following data subject to the
following axioms:

\deflump{Data}
\begin{itemize}
\item
Collection $\ob\bee$ (with elements \emph{0-cells} $A$, $B$,~\ldots)

\item
Categories \trip{A}{B} (with objects \emph{1-cells} $f$, $g$,~\ldots and
arrows \emph{2-cells} $\alpha$, $\beta$,~\ldots)

\item
Functors 
\[
\begin{morphs}
c_{ABC}\from	&\trip{B}{C} \times \trip{A}{B}	&\go	&\trip{A}{C}	\\
		&\pr{g}{f}			&\goesto&g\of f	= gf	\\
		&\pr{\beta}{\alpha}		&\goesto&\beta * \alpha	\\
\end{morphs}
\]
and \gobyc{I_A}{\wun}{\trip{A}{A}} (thus $I_A$ is a 1-cell $A\go A$).

\item
Natural isomorphisms
\[
\begin{ntdiag}
\trip{C}{D}\times\trip{B}{C}\times\trip{A}{B}&	&
\rTo^{1 \times c_{ABC}}	&	&\trip{C}{D}\times\trip{A}{C}	\\
		&	&		&\	&		\\
\dTo<{c_{BCD} \times 1}&&\ruTo>{a_{ABCD}}&	&\dTo>{c_{ACD}}	\\
		&\	&		&	&		\\
\trip{B}{D}\times\trip{A}{B}&&\rTo_{c_{ABD}}&	&\trip{A}{D}	\\
\end{ntdiag}
\]\[
\begin{ntdiag}
\trip{A}{B}\times \wun	&		&	&	&	\\
			&\rdTo(4,4)>{\diso}&	&	&	\\
\dTo<{1\times I_A}	&		&	&	&	\\
			&\ruTo^{r_{AB}}	&	&	&	\\
\trip{A}{B}\times\trip{A}{A}&		&\rTo_{c_{AAB}}&&\trip{A}{B}\\
\end{ntdiag}
\diagspace
\begin{ntdiag}
\wun\times\trip{A}{B} 	&		&	&	&	\\
			&\rdTo(4,4)>{\diso}&	&	&	\\
\dTo<{I_{B}\times 1}	&		&	&	&	\\
			&\ruTo^{l_{AB}}	&	&	&	\\
\trip{B}{B}\times\trip{A}{B}&		&\rTo_{c_{ABB}}&&\trip{A}{B}\\
\end{ntdiag}
\]
thus 2-cells
\[ 
\begin{morphs}
a_{hgf}\from	&(hg)f		&\goiso	&h(gf)	\\
r_f\from	&f \of I_{A}	&\goiso	&f	\\
l_f\from	&I_{B} \of f	&\goiso	&f.	\\
\end{morphs}
\]
\end{itemize}

\deflump{Axioms}
The following commute:
\begin{diagram}
	&	&((kh)g)f&\rTo^{a*1}&(k(hg))f&	&	\\
	&\ldTo<{a}&	&	&	&\rdTo>{a}&	\\
(kh)(gf)&	&	&	&	&	&k((hg)f)\\
	&\rdTo(3,2)<{a}&&	&	&\ldTo(3,2)>{1*a}&\\
	&	&	&k(h(gf))&	&	&	\\
\end{diagram}
\begin{slopeydiag}
(gI)f	&		&\rTo^{a}&		&g(If)	\\
	&\rdTo<{r*1}	&	&\ldTo>{1*l}	&	\\
	&		&gf	&		&	\\
\end{slopeydiag}

\lump{Variants}

If $a$, $l$ and $r$ are identities, so that $(hg)f=h(gf)$, $If=f=fI$, and
similarly for composition of 2-cells, then \bee\ is called a
\emph{2-category}. In this case the axioms hold automatically.

\lump{Example}

There is a 2-category \Cat\ whose 0-cells are small categories, whose 1-cells
are functors, and whose 2-cells are natural transformations.

\lump{Internal Equivalence}

As \Cat\ is a bicategory, we may imitate certain definitions from category
theory in an arbitrary bicategory \bee. In particular, an \emph{(internal)
equivalence} in \bee\ consists of a pair of 1-cells $A \oppair{f}{g} B$
together with an isomorphism $1 \go g \of f$ in the category \trip{A}{A} and
an isomorphism $f \of g \go 1$ in the category \trip{B}{B}. We also say that
$f$ is an equivalence and that $A$ is equivalent to $B$ (inside \bee).

\lump{The Opposite Bicategory}

Given a bicategory \bee, we may form a dual bicategory $\bee^\op$ by
reversing the 1-cells but not the 2-cells. Thus if \bee\ has a 2-cell
\twc{A}{B}{f}{g}{\alpha} then $\bee^\op$ has a 2-cell 
\twcop{A}{B}{f}{g}{\alpha}.

\subsection{Morphisms}

A morphism $F$ (or strictly speaking, \pr{F}{\phi}) from \bee\ to \beep\
consists of the following data subject to the following axioms:

\deflump{Data}
\begin{itemize}
\item
Function \gobyc{F}{\ob\bee}{\ob\beep}

\item
Functors \gobyc{F_{AB}}{\trip{A}{B}}{\tripp{FA}{FB}}

\item
Natural transformations
\[
\begin{ntdiag}
\trip{B}{C}\times\trip{A}{B}&	&\rTo^{c}	&	&\trip{A}{C}	\\
			&	&		&\	&		\\
\dTo<{F_{BC}\times F_{AB}}&	&\ruTo>{\phi_{ABC}}&	&\dTo>{F_{AC}}	\\
			&\	&		&	&		\\
\tripp{FB}{FC}\times\tripp{FA}{FB}&&\rTo_{c'}&	&\tripp{FA}{FC}		\\
\end{ntdiag}
\diagspace
\begin{ntdiag}
\wun	&	&\rTo^{I_A}	&	&\trip{A}{A}	\\
	&	&		&\	&		\\
\dEquals&	&\ruTo>{\phi_{A}}&	&\dTo>{F_{AA}}	\\
	&\	&		&	&		\\
\wun	&	&\rTo_{I'_{FA}}	&	&\tripp{FA}{FA}	\\
\end{ntdiag}
\]
thus 2-cells \gobyc{\phi_{gf}}{Fg \of Ff}{F(g \of f)} and \gobyc{\phi_A}{I'_{FA}}{FI_{A}}.
\end{itemize}

\deflump{Axioms}
The following commute:
\[
\begin{diagram}
(Fh\of Fg)\of Ff	&\rTo^{\phi*1}	&F(h\of g)\of Ff
&\rTo^{\phi}	&F((h\of g)\of f)	\\
\dTo<{a'}		&		&
&		&\dTo>{Fa}		\\
Fh\of (Fg\of Ff)	&\rTo_{1*\phi}	&Fh\of F(g\of f)
&\rTo_{\phi}	&F(h\of(g\of f)) 	\\
\end{diagram}
\]\[
\begin{diagram}
Ff\of I'_{FA}	&\rTo^{1*\phi}	&Ff\of FI_{A}	&\rTo^{\phi}&F(f\of I_A)\\
\dTo<{r'}	&		&		&	&\dTo>{Fr}	\\
Ff		&		&\rEquals	&	&Ff		\\
\end{diagram}
\diagspace
\begin{diagram}
I'_{FB}\of Ff	&\rTo^{\phi*1}	&FI_{B}\of Ff	&\rTo^{\phi}&F(I_B\of f)\\
\dTo<{l'}	&		&		&	&\dTo>{Fl}	\\
Ff		&		&\rEquals	&	&Ff		\\
\end{diagram}
\]

\lump{Variants}

If $\phi_{ABC}$ and $\phi_{A}$ are all natural isomorphisms, so that $Fg\of
Ff \iso F(g\of f)$ and $FI \iso I'$, then $F$ is called a \emph{homomorphism}.
If $\phi_{ABC}$ and $\phi_{A}$ are all identities, so that $Fg\of Ff = F(g\of
f)$ and $FI = I'$, then $F$ is called a \emph{strict homomorphism}.

\lump{Representables}

If $A$ is a 0-cell of a bicategory \bee, there arises a homomorphism
$\gobyc{\rep{A}}{\bee^{\op}}{\Cat}$. The 2-cells ``$\phi_{gf}$'' and
``$\phi_{B}$'' come from $a$ and $r$, respectively.

\lump{Local Properties}

Let $P$ be a property of functors. We say a morphism $F$ is
\emph{locally $P$} if each functor $F_{AB}$ has the property $P$: thus
\emph{locally faithful}, \emph{locally an equivalence},~\ldots.

\subsection{Transformations}

A transformation \twcv{\bee}{\beep}{F}{G}{\sigma}, where $F=\pr{F}{\phi}$ and
$G=\pr{G}{\psi}$ are morphisms, is defined by the following data and axioms.
Below, we use the notation \gobyc{h_*}{\trip{C}{D}}{\trip{C}{E}} for the
functor induced by a 1-cell $D\goby{h}E$ of a bicategory \bee, and similarly
\gobyc{h^*}{\trip{E}{C}}{\trip{D}{C}}.

\deflump{Data}
\begin{itemize}
\item
1-cells $FA \goby{\sigma_A} GA$

\item
Natural transformations
\begin{ntdiag}
\trip{A}{B}	&&\rTo^{F_{AB}}		&&\tripp{FA}{FB}	\\
		&&			&\ &			\\
\dTo<{G_{AB}}	&&\ruTo>{\sigma_{AB}}	&&\dTo>{(\sigma_B)_*}	\\
		&\ &			&&			\\
\tripp{GA}{GB}	&&\rTo_{(\sigma_A)^*}	&&\tripp{FA}{GB}	\\
\end{ntdiag}
thus 2-cells \gobyc{\sigma_f}{Gf\of\sigma_A}{\sigma_{B}\of Ff}.
\end{itemize}

\deflump{Axioms}
The following commute:
\begin{center}
\begin{diagram}
(Gg\of Gf)\of\sigma_{A}	&\rTo^{\wee{a'}}	&
Gg\of(Gf\of\sigma_A)	&\rTo^{1*\sigma_f}	&
Gg\of(\sigma_{B}\of Ff)	&\rTo^{\wee{{a'}^{-1}}}	& 
(Gg\of\sigma_B)\of Ff	&\rTo^{\sigma_{g}*1}	&
(\sigma_{C}\of Fg)\of Ff&\rTo^{\wee{a'}}	&
\sigma_{C}\of(Fg\of Ff)	\\
\dTo<{\psi*1}		&&&&&&&&&& \dTo>{1*\phi}	\\
G(g\of f)\of\sigma_A	&&&&& \rTo_{\sigma_{gf}} &&&&& 
\sigma_{C}\of F(g\of f)	\\
\end{diagram}
\begin{diagram}
I'_{GA}\of\sigma_{A}	&\rTo^{\wee{l'}}&\sigma_A	&
\rTo^{\wee{{r'}^{-1}}}	&\sigma_{A}\of I'_{FA}		\\
\dTo<{\psi*1}		&		&		&	
			&\dTo>{1*\phi}			\\
GI_{A}\of\sigma_A	&		&\rTo_{\sigma_{I_A}}&	
			&\sigma_{A}\of FI_A		\\
\end{diagram}
\end{center}

\lump{Variants}

If $\sigma_{AB}$ are all natural isomorphisms then $\sigma$ is called a
\emph{strong transformation}. If $\sigma_{AB}$ are all identities then
$\sigma$ is called a \emph{strict transformation}.

\lump{Representables}

If $A \goby{f} B$ is a 1-cell in a bicategory, then there arises a
strong transformation \gobyc{f_{*}=\rep{f}}{\rep{A}}{\rep{B}}. The 2-cells
``$\sigma_g$'' come from the associativity isomorphism $a$.

\subsection{Modifications}

A modification 
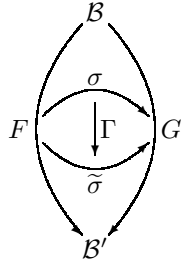
\[
\thcv{\bee}{\beep}{F}{G}{\sigma}{\twid{\sigma}}{\Gamma}
\]
consists of the following data subject to the following axioms:

\deflump{Data}
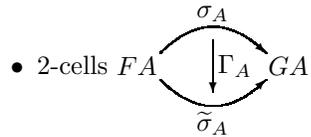
\begin{itemize}
\item
2-cells \twc{FA}{GA}{\sigma_A}{\twid{\sigma}_A}{\Gamma_A}
\end{itemize}

\deflump{Axioms}
The following commute:
\begin{diagram}
Gf\of\sigma_A	&\rTo^{1*\Gamma_A}	&Gf\of\twid{\sigma}_A	\\
\dTo<{\sigma_f}	&			&\dTo>{\twid{\sigma}_f}	\\
\sigma_{B}\of Ff&\rTo_{\Gamma_{B}*1}	&\twid{\sigma}_{B}\of Ff\\
\end{diagram}

\lump{Variants}

None.

\lump{Representables}

If \twc{A}{B}{f}{g}{\alpha} is a 2-cell in a bicategory, then there arises a
modification 
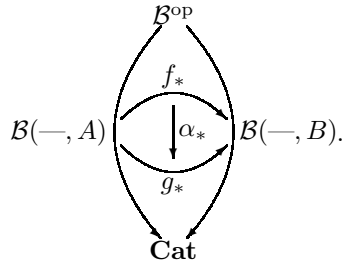
\[
\thcv{\bee^{\op}}{\Cat}{\rep{A}}{\rep{B}}{f_*}{g_*}{\alpha_*}.
\]

\subsection{Strength Terminology}

The terminology used to describe whether something holds strictly, up to
isomorphism, or just up to a connecting map, has evolved messily.  Here is a
summary of the definitions given above; note that for representables,
everything is at the `iso' level.
\begin{center}
\begin{tabular}{c|ccc}
	&bicategories	&morphisms		&transformations	\\
\hline
map	&---		&(plain)		&(plain)		\\
iso	&bicategory	&homomorphism		&strong	transformation	\\
equality&2-category	&strict homomorphism	&strict transformation
\end{tabular}
\end{center}



\section{Coherence}
\setcounter{subsection}{-1}

\subsection{Functor Bicategories}

Given a pair of bicategories \bee\ and \beep, one can define (canonically) a
`functor bicategory' \homset{\mb{Lax}}{\bee}{\beep}, whose 0-cells are
morphisms \bee\go\beep, whose 1-cells are transformations, and whose
2-cells are modifications. This is not in general a 2-category, but it is if
\beep\ is. We will take a particular interest in the sub-bicategory
\ftrbi{\bee}{\beep} of \homset{\mb{Lax}}{\bee}{\beep}, consisting of
homomorphisms, strong transformations and modifications (i.e.\ everything at
the `iso' level). 

\subsection{The Yoneda Embedding}	\label{Yoneda}

By the observations of the previous section, there is for each bicategory
\bee\ a 2-category \Pshf. Moreover, the representables constructed in
Section~\ref{sec:defns} provide a `Yoneda' homomorphism
\gobyc{Y}{\bee}{\Pshf}. It is straightforward to calculate that locally $Y$
is full, faithful, and essentially surjective on objects---in other words,
that $Y$ is a local equivalence. (In order for $Y$ to exist it is necessary
that \bee\ should be locally small, but we do not emphasize this issue here.)

\subsection{Biequivalence}

Let \bee\ and \beep\ be bicategories. A \emph{biequivalence} from \bee\ to
\beep\ consists of a pair of homomorphisms $\bee\oppair{F}{G}\beep$
together with an equivalence $1 \go G \of F$ inside the bicategory 
\ftrbi{\bee}{\bee} and an equivalence $F \of G \go 1$ inside
\ftrbi{\beep}{\beep}. We also say that $F$ is a biequivalence and that \bee\
is biequivalent to \beep. Now, just as for equivalence of plain categories,
there is an alternative criterion for biequivalence: namely, that a
homomorphism \gobyc{F}{\bee}{\beep} is a biequivalence if and only if $F$ is
locally an equivalence and is surjective-up-to-equivalence on objects. The
latter condition means that if $B'$ is any 0-cell of \beep\ then there is
some 0-cell $B$ of \bee\ such that $FB$ is (internally) equivalent to $B'$.

\subsection{The Coherence Theorem}	\label{CohThm}

{\raggedright
\textbf{Theorem}
\textit{Every bicategory is biequivalent to a 2-category.}}
\vspace{1ex} \newline
\pf\ Let \bee\ be a bicategory and \gobyc{Y}{\bee}{\Pshf} the Yoneda map. Let
\beep\ be the full image of $Y$: that is, the sub-2-category of \Pshf\
whose 0-cells are those in the image of $Y$, and with all 1- and 2-cells of
\Pshf\ between them. Let \gobyc{Y'}{\bee}{\beep} be the restriction of $Y$. 
Then:
\begin{itemize}
\item
$Y'$ is a homomorphism, since $Y$ is
\item
$Y'$ is surjective on 0-cells, by construction
\item
$Y'$ is a local equivalence, since $Y$ is (observed in~\ref{Yoneda}).
\end{itemize}
Thus $Y'$ is a biequivalence from \bee\ to the 2-category \beep.
\done

\subsection{Coherence and Commuting Diagrams}
\label{sec:coh-and-comm}

Coherence theorems sometimes have the form `all diagrams of a certain kind
commute'; for instance, the coherence theorem for monoidal categories states
that all diagrams built out of the associativity and identity isomorphisms
commute. This result for monoidal categories also holds for bicategories.
Without giving a precise statement or proof, we indicate by an example how it
is a corollary of our coherence theorem,~\ref{CohThm}.

A typical instance of what we wish to prove is that the diagram
\begin{equation}	\label{above}
\begin{diagram}
		&		&(h(Ig))f&		&	\\
		&\ldTo<{a^{-1}*1}&	&\rdTo(2,3)>{a}	&	\\
((hI)g)f	&		&	&		&	\\
\dTo<{(r*1)*1}	&		&	&		&h((Ig)f)\\
(hg)f		&		&	&\ldTo(2,3)>{1*(l*1)}&	\\
		&\rdTo<{a}	&	&		&	\\
		&		&h(gf)	&		&	\\
\end{diagram}
\end{equation}
should commute, for any composable 1-cells $f$, $g$, $h$ in any bicategory.
Let us say that a bicategory \bee\ has the \emph{coherence property} if all
diagrams `like this' in \bee\ commute; our goal is to show that every
bicategory has the coherence property. To achieve this, first observe that
every 2-category has the coherence property, since $a$, $l$ and $r$ are all
1. Then, for any \bee\ we have an `embedding' \gobyc{Y}{\bee}{\Pshf} of
\bee\ into a bicategory with the coherence property, and this implies that
\bee\ too has the coherence property, as now explained.

Let \gobyc{\pr{F}{\phi}}{\bee}{\beep} be a morphism of bicategories, and
suppose that \beep\ has the coherence property. We want to deduce that,
subject to certain conditions on \pr{F}{\phi}, the bicategory \bee\ also has
the coherence property. Consider the diagram~(\ref{above}) in \bee. Let
$\alpha$ be the composite down the left-hand side, and $\beta$ down the
right; let $\alpha'$ be the composite
\[
(Fh\of(I'\of Fg))\of Ff \goby{{a'}^{-1}*1}
((Fh\of I')\of Fg)\of Ff \goby{(r'*1)*1}
(Fh\of Fg)\of Ff \goby{a'} Fh\of (Fg\of Ff),
\]
and similarly $\beta'$. We know that $\alpha'=\beta'$ and want to conclude
that $\alpha=\beta$. Consider, then, the diagram in Figure~2,
\setcounter{figure}{1}
\begin{figure}
\begin{diagram}
(Fh\of(I'\of Fg))\of Ff	&\rTo^{\phi}	&F[(h(Ig))f]		\\
\dTo<{{a'}^{-1}*1}	&		&\dTo>{F[a^{-1}*1]}	\\
((Fh\of I')\of Fg)\of Ff&\rTo^{\phi}	&F[((hI)g)f]		\\
\dTo<{(r'*1)*1}		&		&\dTo>{F[(r*1)*1]}	\\
(Fh\of Fg)\of Ff	&\rTo^{\phi}	&F[(hg)f]		\\
\dTo<{a'}		&		&\dTo>{F[a]}		\\
Fh\of (Fg\of Ff)	&\rTo_{\phi}	&F[h(gf)]		\\
\end{diagram}
\caption{diagram of 2-cells in \beep}
\end{figure}
where the 2-cells called $\phi$ are built up from $\phi_A$'s and
$\phi_{qp}$'s. By definition of morphism, the diagram commutes---that is,
$F\alpha\of\phi = \phi\of\alpha'$. Similarly, $F\beta\of\phi =
\phi\of\beta'$, so $F\alpha\of\phi = F\beta\of\phi$. If $F$ is a homomorphism
then the 2-cells $\phi$ are isomorphisms, so $F\alpha = F\beta$; if also $F$
is locally faithful then we may conclude that $\alpha=\beta$.

Our example therefore demonstrates: if \gobyc{F}{\bee}{\beep} is a locally
faithful homomorphism and \beep\ has the coherence property, then so does
\bee. Applying this to \gobyc{Y}{\bee}{\Pshf}, for any bicategory \bee, shows
that all bicategories have the coherence property.





\begin{thebibliography}{M}	\raggedright
\addcontentsline{toc}{section}{\numberline{}References}


\bibitem{Ben}
Jean B\'{e}nabou,
Introduction to bicategories (1967).
In: \emph{Reports of the Midwest Category Seminar}, ed.\ B\'{e}nabou et
al, Springer LNM 47. 

\bibitem{GPS}
R. Gordon, A. J. Power, Ross Street,
Coherence for tricategories (1995).
\emph{Memoirs of the AMS}, vol.~117, no.~558.

\bibitem{Gray}
J. Gray,
\emph{Formal Category Theory: Adjointness for 2-Categories}\/ (1974).
Springer LNM 391.

\bibitem{KS}			
G. M. Kelly, R. Street,
Review of the elements of 2-categories (1974).
In: \emph{Category Seminar}, Springer LNM 420, pp~75--103.

\bibitem{Lack}
Stephen Lack,
\emph{The Algebra of Distributive and Extensive Categories} (1995).
PhD thesis, University of Cambridge. Available via
\url{http://cat.maths.usyd.edu.au/$\sim$stevel}. 

\bibitem{StFB}
Ross Street,
Fibrations in bicategories (1980).
\emph{Cahiers Top.\ Geom.\ Diff.}, vol.~XXI, no.~2. 

\bibitem{StCS}
Ross Street,
Categorical structures (1995).
In: \emph{Handbook of Algebra,}\/ Vol.~1, ed. M. Hazewinkel, Elsevier
North-Holland. 


\end{thebibliography}
\end{document}